\documentclass[eng]{isaac}
\usepackage{amsmath, amstext, amsfonts}

\begin{document}
\title{COGARCH: Symbol, Generator and Characteristics}

% Authors
% First author block
\author{A. Schnurr}
\address{TU Dortmund, Germany, 44227 Dortmund, Fakult\"at f\"ur Mathematik, Vogelpothsweg 87, tel: +49-231-755-3099, e-mail: alexander.schnurr@math.tu-dortmund.de  }
% Second author block
%\author{A. N. SecondAuthor}
%\address{Contacts for the second author: Company, Country, ZIP, City, address, phone, E-mail}
%


\amsnum{60J75 (primary), 60G51, 60J25 (secondary)}

\abstracts{ We describe the technique how to use the symbol in order to calculate the generator and the characteristics of an It\^o process. As an example we analyze the COGARCH process which is used to model financial data.  }

\maketitle

\section{Introduction}

%\texttt{amsmath},\texttt{amstext} \texttt{amsfonts}. http://www.ams.org/tex/amslatex.html+.

The COGARCH process was introduced by Kl\"uppelberg et al. in \cite{cogarch} in order to model financial data. It is a continuous time analog of the classic GARCH process (in discrete time) and it is based on a single background driving L\'evy process in contrast to the well known model by Barndorff-Nielsen and Shephard \cite{bar-she01}. L\'evy processes are c\`adl\`ag universal Markov processes which are homogeneous in time \emph{and} space. Our main reference for this class of processes is \cite{sato}. For the L\'evy triplet we write ($\ell, Q, N$).

In the present paper we calculate the so called \emph{symbol} of the COGARCH process (and its volatility process). The origins of the symbols are in the theory of partial differential equations, namely they appear in the Fourier representation of certain operators. The symbol found its way into probability theory for the following reason: suppose we are given a Feller process $X$ with associated semigroup $(T_t)_{t\geq 0}$ and generator $(A,D(A))$. Suppose further that the test functions $C_c^\infty(\Bbb{R}^d)$ are contained in the domain $D(A)$. In this case $A$ is a pseudo-differential operator with symbol $-q(x,\xi)$. For every $x\in\Bbb{R}^d$ $q(x,\cdot)$ is a continuous negative definite function in the sense of Schoenberg (cf. \cite{bergforst} Chapter 2).

For a detailed, self contained treatment on the interplay between the process and its symbol cf. the monograph \cite{niels}. In this context the following four questions are of interest: \\
I) Given a process, (say as the solution of an SDE) what is its symbol? (E.g. \cite{SDEsymbol})\\
II) Given a symbol, does there exist a corresponding process? (\cite{hoh00, Hoh2002, jacobschilling})\\
III) Which properties of the process can be characterized via the symbol? ( \cite{schilling98hdd, schilling98})\\
IV) For which bigger classes of processes is it possible (and useful) to define a symbol? (\cite{generalsymbindices, mydiss}) \\
All four questions are a vital part of ongoing research. In the present paper we emphasize, how one can calculate the symbol of a given process using a probabilistic formula and derive directly the generator as well as the semimartingale characteristics. 

The notation we are using is (more or less) standard. Vectors are meant to be column vectors and the transposed of a vector $v$ or a matrix $Q$ is denoted by $v'$ respective $Q'$. 

Let us recall how the COGARCH process is defined:\newline
we start with a L\'evy process $Z=(Z_t)_t$ with triplet $(\ell,Q,N)$. Fix $0<\delta< 1, \ \beta>0,\ \lambda \geq 0$.
Then the volatility process $(\sigma_t)_{t\geq0}$ is the solution of the SDE
\begin{eqnarray*}
d\sigma^2_t &=& \beta \ dt + \sigma_t^2 \left( \log \delta \ dt + \frac{\lambda}{\delta} \ d[Z,Z]_t^{disc} \right)\\
\sigma_0    &=&S
\end{eqnarray*}
where $S>0$ and
\begin{eqnarray*}
[Z,Z]_t^{disc}=\sum_{0<s\leq t} (\Delta Z_s)^2.
\end{eqnarray*}
It turns out, that $(\sigma_t)_{t\geq 0}$ is a time homogeneous Markov process.

\textbf{Definition:}
The process
\begin{eqnarray*}
G_t:=g+\int_0^t \sigma_{s-} \ dZ_t, \hspace{1cm} g\in \Bbb{R},
\end{eqnarray*}
is called \textbf{COGARCH process} (starting in $g$). 

We allow the process to start everywhere in order to bring our methods into account. The pair $(G_t, \sigma_t^2)$ is a (normal) Markov process which is homogeneous in time. It is homogeneous in space in the first component. 
Furthermore $(G_t, \sigma_t^2)$ is an It\^o process, which follows from Theorem 3.33 of \cite{cinlarjacod81} which characterizes It\^o processes as solutions of certain stochastic differential equations and Proposition IX.5.2. of \cite{jacodshir} giving a representation of the semimartingale characteristics of a stochastic integral.

To avoid problems which might arise for processes defined on $\Bbb{R} \times \Bbb{R}_+$ we consider in the following: $(G_t,V_t)=(G_t, \log(\sigma_t^2))$, i.e. $V$ is the logarithmic squared volatility. 

%%% symbol %%%%%%%%%%%%%%%%%%%%%%%%%%%%%%%%%%%%%%%%%%%%%%%%%%%%%%%%%%%%%%%%%%%%%%%%%%%%
\section{The Symbol of a Stochastic Process}

\textbf{Definition:}
Let $X$ be an $\Bbb{R}^d$-valued universal Markov process, which is conservative and normal. Fix a starting point $x$ and define $T=T^x_R$ to be the first exit time from the ball of radius $R > 0$: 
  \begin{eqnarray} \label{stoppingtime}
     T:=T^x_R:=\inf\{t\geq 0 : \left\|X_t-x\right\| > R \} \text{ under } \Bbb{P}^x (x\in\Bbb{R}^d).
  \end{eqnarray}
We call the function $p:\Bbb{R}^d\times \Bbb{R}^d \rightarrow \Bbb{C}$, given by
\begin{eqnarray} \label{symbol}
   p(x,\xi):=- \lim_{t\downarrow 0}\Bbb{E}^x \frac{e^{i(X^T_t-x)'\xi}-1}{t},
\end{eqnarray}
the \textbf{(probabilistic) symbol} of the process, if the limit exists for every $x,\xi$ and $R$ and is independent of the choice of $R$. 

In \cite{mydiss} Theorem 4.4. we have shown that for It\^o processes in the sense of Cinlar, Jacod, Protter and Sharpe (cf. \cite{vierleute}) having differential characteristics which are finely continuous (cf. \cite{blumenthalget}) and locally bounded the above limit exists and coincides for every choice of $R$. For the reader's convenience we recall the the definition of It\^o processes, as it is used here:

\textbf{Definition:}
A Markov semimartingale $X=(X_t)_{t\geq 0}$, i.e. a universal Markov process which is a semimartingale w.r.t. every initial probability $\Bbb{P}^x$ ($x\in\Bbb{R})$, is called \textbf{It\^o process} if it has characteristics of the form:
\begin{eqnarray*}
B_t^{j}(\omega)=&\int_0^t \ell^{j}(X_s(\omega)) \ ds & j=1,...,d\\
C_t^{jk} (\omega)=&\int_0^t Q^{jk} (X_s(\omega)) \ ds &j,k=1,...,d\\
\nu(\omega;ds,dy)=&N(X_s(\omega),dy)\ ds 
\end{eqnarray*}
where $\ell^{j},Q^{jk}:\Bbb{R}^d \to \Bbb{R}$ are measurable functions, $Q(x)=(Q^{jk}(x))_{1\leq j,k \leq d}$ is a positive semidefinite matrix for every $x\in\Bbb{R}^d$, and $N(x,\cdot)$ is a Borel transition kernel on $\Bbb{R}^d \times \mathcal{B}(\Bbb{R}^d \backslash \{0\})$. $\ell$, $Q$ and $\int_{y\neq 0}(1\wedge y^2)N(\cdot,dy)$ are called \textbf{differential characteristics}.

Example 1:
Let $X$ be a $d$-dimensional L\'evy process. It is a well known fact that the characteristic function of $X_t$ ($t\geq 0$) can be written as
\begin{eqnarray*}
\Bbb{E}^0 \exp(iX_t'\xi)= \exp(-t\psi(\xi)).
\end{eqnarray*}
The function $\psi:\Bbb{R}^d\to \Bbb{C}$ is called characteristic exponent. By an elementary calculation one obtains $p(x,\cdot)=\psi(\cdot)$ for every $x\in\Bbb{R}^d$.

Example 2:
Let $X$ be a rich Feller process, i.e. the test functions $C_c^\infty(\Bbb{R}^d)$ are contained in the domain $D(A)$ of the generator $A$. In this case the generator restricted to $C_c^\infty(\Bbb{R}^d)$ is a pseudo-differential operator with (functional analytic) symbol $-q(x,\xi)$. In \cite{mydiss} we have shown that $X$ is an It\^o process and $p(x,\xi)=q(x,\xi)$ for every $x,\xi\in\Bbb{R}^d$.

Example 3:
Let $(Z_t)_{t\geq 0}$ be an $\Bbb{R}^n$-valued L\'evy process. The solution of the stochastic differential equation ($x\in\Bbb{R}^d$),
\begin{eqnarray*}
  dX_t^x&=&\Phi(X_{t-}^x) \,dZ_t \\
  X_0^x&=&x, 
\end{eqnarray*}
where $\Phi: \Bbb{R}^d \to \Bbb{R}^{d \times n}$ is Lipschitz continuous admits the symbol
\begin{eqnarray*}
    p(x,\xi)=\psi(\Phi(x)'\xi).
\end{eqnarray*}
This was shown in \cite{SDEsymbol}.

%%% generator, char %%%%%%%%%%%%%%%%%%%%%%%%%%%%%%%%%%%%%%%%%%%%%%%%%%%%%%%%%%%%%%%%%%%
\section{Symbol, Generator and Characteristics}

In the present section we calculate the symbol of the COGARCH process. Using the close relationship between the symbol, the extended generator and the semimartingale characteristics we are able to write down the latter two objects directly. Let us emphasize that the symbol does \emph{not} depend on $g$, since the process is homogeneous in the first component. 

\textbf{Theorem:}
The stochastic process $(G_t,V_t)=(G_t, \log(\sigma_t^2))$ admits the symbol $p:\Bbb{R}^2 \times \Bbb{R}^2 \to \Bbb{C}$ given by \footnotesize
\begin{eqnarray*}
&&p\left( \binom{g}{v} ,\xi \right) =\\
&&\hspace{2mm} -i \xi_1 \left(\ell e^{v/2} + e^{v/2} \int_{\Bbb{R} \backslash \{ 0\}} y \cdot (1_{\{\left|e^{v/2}y\right|<1\}} \cdot 1_{\{ \left|\log(1+(\lambda/\delta) \ y^2)\right|<1 \}} - 1_{\{|y| < 1\}}) \ N(dy) \right)\\
&&\hspace{2mm}- i\xi_2 \left(\frac{\beta}{e^{v}} + \log \delta + \int_{\Bbb{R} \backslash \{ 0\}} \log(1+\frac{\lambda}{\delta}y^2) \cdot (1_{\{\left|e^{v/2}y\right|<1\}} \cdot 1_{\{ \left|\log(1+(\lambda/\delta) \ y^2)\right|<1 \}} ) \ N(dy) \right) \\
&&\hspace{2mm}+\frac{1}{2} \xi_1^2 e^{v} Q \\
&&\hspace{2mm}-\int_{\Bbb{R}^2 \backslash \{0 \}} \Big( e^{i(z_1,z_2)\xi}-1-i z'\xi \cdot (1_{\{\left|z_1\right|<1\}} \cdot 1_{\{ \left|z_2\right|<1 \}} ) \Big) 
  \tilde{N} \left( \binom{g}{v} , dz \right), \\
\end{eqnarray*} \normalsize
where $\tilde{N}$ is the image measure
\begin{eqnarray*}
\tilde{N} \left( \binom{g}{v} , dz \right) = N(f_v\in dz)
\end{eqnarray*}
under $f:\Bbb{R}\to\Bbb{R}^2$ given by
\begin{eqnarray*}
f_v(w)= \binom{e^{v/2}w}{\log(1+(\lambda/\delta) \ w^2)}.
\end{eqnarray*}

Remark: It is not surprising, that the transformation of the jump measure depends only on $v$ since the process is space homogeneous in the first component.

\textbf{Proof:} Let $T$ be the stopping time defined in \eqref{stoppingtime}. At first we use It\^o's formula:
\begin{align}
& \frac{\Bbb{E}^{g,v}e^{i (G_t^T-g,V_t^T-v)\xi}-1}{t} 
= \frac{\Bbb{E}^{0,v} e^{i (G_t^T,V_t^T-v)\xi} -1}{t}  \nonumber \\
&= \frac{1}{t} \Bbb{E}^{0,v} \int_{0+}^t i\xi_1  e^{i (G_{s-}^T,V_{s-}^T-v)\xi} \ dG_s^T \tag{I}\\
&+ \frac{1}{t} \Bbb{E}^{0,v} \int_{0+}^t i\xi_2  e^{i (G_{s-}^T,V_{s-}^T-v)\xi} \ dV_s^T \tag{II}\\
&- \frac{1}{2t} \Bbb{E}^{0,v} \int_{0+}^t \xi_1^2  e^{i (G_{s-}^T,V_{s-}^T-v)\xi} \ d[G^T,G^T]_s^c\tag{III} \\
&- \frac{1}{t} \Bbb{E}^{0,v} \int_{0+}^t \xi_1\xi_2  e^{i (G_{s-}^T,V_{s-}^T-v)\xi} \ d[G^T,V^T]_s^c \tag{IV}\\
&- \frac{1}{2t} \Bbb{E}^{0,v} \int_{0+}^t \xi_2^2  e^{i (G_{s-}^T,V_{s-}^T-v)\xi} \ d[V^T,V^T]_s^c \tag{V}\\
&+  \frac{1}{t}\Bbb{E}^{0,v}\sum_{0<s\leq t} e^{(G_{s-}^T,V_{s-}^T-v)\xi} 
\left(e^{i\Delta(G_s^T,V_s^T)\xi} -1 -(i\xi_1 \Delta G_s^T + i\xi_2 \Delta V_s^T)\right). \tag{VI}
\end{align}
We deal with this formula term-by-term. In the calculation of the first term we use
\begin{eqnarray*}
dG_s^T=\sigma_{s-} 1_{\{s\in\left[ \left[ 0,T \right] \right]\}} \ dZ_s.
\end{eqnarray*}
Recall that the integrand is bounded and for the L\'evy process $Z$ we have the L\'evy-It\^o-decomposition:
\begin{align*}
    Z_t=\ell t&+\sqrt{Q} W_t
        +\int_{[0,t] \times \{\left|y\right|<1\}} y \ (\mu^Z(ds,dy)-ds N(dy)) \\
        &+\sum_{0<s\leq t} \Delta Z_s 1_{\{\left|\Delta Z_s\right| \geq 1 \}  },
\end{align*}
where $\mu^Z$ denotes the jump measure of the process (cf. \cite{jacodshir} Proposition II.1.16).
The integrals with respect to the martingale parts are again $L^2$-martingales and the respective terms disappear. What remains from the first term is:
\begin{eqnarray} \label{termone}
\frac{1}{t} \Bbb{E}^{0,v} \int_{0+}^t i\xi_1  e^{i (G_{s-}^T,V_{s-}^T-v)\xi} \sigma_{s-} 1_{\{s\in\left[ \left[ 0,T \right] \right]\}} \ d\left(\ell s + \sum_{0<r\leq s} \Delta Z_r \cdot 1_{\{\left|\Delta Z_r\right|\geq 1\}}\right).
\end{eqnarray}
For the first part of this integrand we get:
\begin{eqnarray*}
&&\frac{1}{t} \Bbb{E}^{0,v} \int_{0+}^t i\xi_1  e^{i (G_{s-}^T,V_{s-}^T-v)\xi} \sigma_{s-} 1_{\{s\in\left[ \left[ 0,T \right] \right]\}} \ d(\ell s) \\
&&\hspace{10mm}= \Bbb{E}^{0,v} \frac{1}{t}\int_{0}^t i\xi_1 \ell e^{i (G_{s}^T,V_{s}^T-v)\xi}  1_{\{s\in\left[ \left[ 0,T \right[ \right[\}}\sigma_{s} \ ds \\
&&\hspace{10mm}= i\xi_1 \ell \ \Bbb{E}^{0,v} \underbrace{\int_{0}^1  e^{i (G_{st}^T,V_{st}^T-v)\xi}  1_{\{st\in\left[ \left[ 0,T \right[ \right[\}}}_{\to 1} \underbrace{\sigma_{st}}_{\to S} \ ds \\
&&\xrightarrow[t\downarrow 0]{} i \xi_1 \ell S.
\end{eqnarray*}
In the first equation we used the fact that we are integrating with respect to Lebesgue measure. For this the countable number of jump times is a nullset. In the last step we used Lebesgue's theorem twice. A similar argumentation is used in the consideration of the second and the third term. The jump term of (\ref{termone}) above will be compared to the sixth term. 

Using It\^o's formula  we obtain for the second term \footnotesize
\begin{eqnarray*}
\frac{1}{t} \Bbb{E}^{0,v} \int_{0+}^t i\xi_2  e^{i (G_{s-}^T,V_{s-}^T-v)\xi} \left\{\frac{1}{\sigma_{s-}^2} \ d(\sigma_s^T)^2 + d\left(\sum_{0<r\leq s} \log \sigma_r^2 - \log \sigma_{r-}^2 - \frac{1}{\sigma_{r-}^2} \Delta(\sigma_r^2) \right) \right\}
\end{eqnarray*} \normalsize
and by plugging in the defining SDE for $(\sigma^2)$:
\begin{eqnarray*}
&&\frac{1}{t} \Bbb{E}^{0,v} \int_{0+}^t i\xi_2  e^{i (G_{s-}^T,V_{s-}^T-v)\xi}1_{\{s\in\left[ \left[ 0,T \right] \right]\}} \ 
\left\{ \rule[0mm]{0mm}{7mm}\left( \frac{\beta}{\sigma_{s-}^2}  \ ds + \frac{\sigma_{s-}^2}{\sigma_{s-}^2}  \log \delta \ ds \right) \right.\\
&&\left.+ \frac{\lambda}{\delta} \ d\left(\sum_{0<r\leq s} (\Delta Z_r)^2 )\right) +
d\left(\sum_{0<r\leq s} \Delta(\log \sigma_r^2) - \frac{1}{\sigma_{r-}^2} \Delta(\sigma_r^2) \right) \right\}.
\end{eqnarray*}
We postpone the jump parts and for the remainder term we get in the limit, using a similar argumentation as for the first term,
\begin{eqnarray*}
\xrightarrow[t\downarrow 0]{}i\xi_2 \beta / S^2 + i \xi_2 \log \delta.
\end{eqnarray*}
For the third term we obtain in an analogous manner to the first one
\begin{eqnarray*}
&&-\frac{1}{2t} \Bbb{E}^{0,v} \int_{0+}^t \xi_1^2  e^{i (G_{s-}^T,V_{s-}^T-v)\xi} \ d[G^T,G^T]_s^c \\
&&\hspace{1cm}= -\frac{1}{2t} \Bbb{E}^{0,v} \int_{0+}^t \xi_1^2  e^{i (G_{s-}^T,V_{s-}^T-v)\xi} \  1_{\{s\in\left[ \left[ 0,T \right] \right]\}} \sigma_{s-}^2 \ d[Z,Z]_s^c \\
&&\hspace{1cm}= -\frac{1}{2t} \Bbb{E}^{0,v} \int_{0}^t \xi_1^2  e^{i (G_{s-}^T,V_{s-}^T-v)\xi} \  1_{\{s\in\left[ \left[ 0,T \right[ \right[\}} \sigma_{s-}^2 \ d(Qs) \\
&&\xrightarrow[t\downarrow 0]{} -\frac{1}{2} \xi_1^2 S^2 Q.
\end{eqnarray*}
The terms four and five are constant zero: since $(t)_t$ and $([Z,Z]_t)_t$ are both of finite variation on compacts, the process $(\sigma_t^2)_t$ has this property as well, by its very definition. Therefore it is a quadratic pure jump process (see \cite{protter} Section II.6). Using It\^o's formula we obtain that $V=\log(\sigma^2)$ is again a quadratic pure jump process and therefore
\begin{eqnarray*}
[V^T, V^T]^c_s=0 \text{ and } [V^T, G^T]^c_s=0 .
\end{eqnarray*}
The only thing that remains to do is dealing with the various `jump parts'.
From the first term we left the following behind
\begin{eqnarray*}
&&\frac{1}{t} \Bbb{E}^{0,v} \int_{0+}^t i\xi_1  e^{i (G_{s-}^T,V_{s-}^T-v)\xi} \sigma_{s-} 1_{\{s\in\left[ \left[ 0,T \right] \right]\}} \ d\left( \sum_{0<r\leq s} \Delta Z_r \cdot 1_{\{\left|\Delta Z_r\right|\geq 1\}}\right) \\
&&\hspace{1cm}= \frac{1}{t} \Bbb{E}^{0,v}  \sum_{0<s\leq t}i\xi_1  e^{i (G_{s-}^T,V_{s-}^T-v)\xi} \sigma_{s-} 1_{\{s\in\left[ \left[ 0,T \right] \right]\}}  \Delta Z_s \cdot 1_{\{\left|\Delta Z_s\right|\geq 1\}}
\end{eqnarray*}
and from the second one
\begin{eqnarray*}
&&\frac{1}{t} \Bbb{E}^{0,v} \int_{0+}^t i\xi_2  e^{i (G_{s-}^T,V_{s-}^T-v)\xi}1_{\{s\in\left[ \left[ 0,T \right] \right]\}} \frac{\lambda}{\delta} \ d\left(\sum_{0<r\leq s} (\Delta Z_r)^2\right) 
\\
&&  \hspace{1cm}+ \frac{1}{t} \Bbb{E}^{0,v} \int_{0+}^t i\xi_2  e^{i (G_{s-}^T,V_{s-}^T-v)\xi}1_{\{s\in\left[ \left[ 0,T \right] \right]\}} \
d\left(\sum_{0<r\leq s} \Delta V_r - \frac{1}{\sigma_{r-}^2} \Delta(\sigma_r^2) \right) \\
&&=\frac{1}{t} \Bbb{E}^{0,v} \sum_{0<s\leq t} i\xi_2  e^{i (G_{s-}^T,V_{s-}^T-v)\xi}1_{\{s\in\left[ \left[ 0,T \right] \right]\}} \frac{\lambda}{\delta} (\Delta Z_s)^2
\\
&&  \hspace{1cm}+ \frac{1}{t} \Bbb{E}^{0,v} \sum_{0<s\leq t} i\xi_2  e^{i (G_{s-}^T,V_{s-}^T-v)\xi}1_{\{s\in\left[ \left[ 0,T \right] \right]\}} 
 \left( \Delta V_s - \frac{1}{\sigma_{s-}^2} \Delta(\sigma_s^2) \right).
\end{eqnarray*}
Adding these terms to term number six and using the equalities
\begin{eqnarray*}
\Delta G_s^T = (\sigma_{s-} 1_{\{s\in\left[ \left[ 0,T \right] \right]\}}) \Delta Z_s \text{ and }
(\Delta \sigma_s^T)^2 = \frac{\lambda}{\delta}(\sigma_{s-}^2 1_{\{s\in\left[ \left[ 0,T \right] \right]\}}) (\Delta Z_s)^2
\end{eqnarray*}
as well as
\begin{eqnarray*}
\Delta \log(\sigma_s^2)^T = \log \left( \frac{(\sigma_{s-}^2)^T + \Delta(\sigma_{s}^2)^T}{(\sigma_{s-}^2)^T} \right) = \log\left(1+\frac{\Delta(\sigma_{s}^2)^T}{(\sigma_{s-}^2)^T} \right)
\end{eqnarray*}
we obtain
\begin{eqnarray*} \footnotesize
&&\frac{1}{t} \Bbb{E}^{0,v}  \sum_{0<s\leq t}  e^{i (G_{s-}^T,V_{s-}^T-v)\xi}  
  1_{\{s\in\left[ \left[ 0,T \right] \right]\}} \ \times\\
&&\hspace{1cm}\left( e^{i\sigma_{s-} \Delta Z_s \xi_1 + i\log(1+ (\lambda / \delta)\Delta (Z_s)^2) \xi_2} -1 - 
  i\xi_1 \sigma_{s-} \Delta Z_s \cdot 1_{\{\left|\Delta Z_s\right| < 1\}} \right) \\
&&=\frac{1}{t} \Bbb{E}^{0,v}  \int_{\left] 0,t \right]\times \{y\neq 0\}}e^{i 
  (G_{s-}^T,V_{s-}^T-v)\xi}  1_{\{s\in\left[ \left[ 0,T \right] \right]\}} \ \times \\
&&\hspace{1cm}\left( e^{i\sigma_{s-} y \xi_1 + i\log(1+ (\lambda / \delta)y^2) \xi_2} -1 - i\xi_1 \sigma_{s-} 
  y \cdot 1_{\{|y| < 1\}} \right)\ \mu^Z(\cdot;ds,dy) \\
&&=\frac{1}{t} \Bbb{E}^{0,v}  \int_{\left] 0,t \right]\times \{y\neq 0\}}e^{i   
  (G_{s-}^T,V_{s-}^T-v)\xi}  1_{\{s\in\left[ \left[ 0,T \right] \right]\}} \ \times \\
&&\hspace{-7mm}\left(\rule[8mm]{0mm}{0mm}\left( e^{i\sigma_{s-} y \xi_1 + i\log(1+ (\lambda / \delta)y^2) \xi_2} -1 - i\binom{\sigma_{s-}y}{ \log(1+\frac{\lambda}{\delta} y^2)}'\xi\cdot 1_{\{|Sy|<1\}} \cdot 1_{\{ 
  \left|\log(1+\frac{\lambda}{\delta} y^2)\right|<1 \}}  \right)  \right.\\
&&\hspace{1cm}+\left( i\xi_1 \sigma_{s-}y \cdot (1_{\{|Sy|<1\}} \cdot 1_{\{ \left|\log(1+\frac{\lambda}{\delta} y^2)\right| <1\}} ) 
  - 1_{\{|y| < 1\}}) \right) \\
&&\hspace{1cm} \left. \rule[8mm]{0mm}{0mm} +\left(i\xi_2 \log(1+\frac{\lambda}{\delta} y^2) \cdot 1_{\{|Sy|<1\}} \cdot 1_{\{ \left|\log(1+\frac{\lambda}{\delta} y^2)\right|<1 \}} \right) \right) \ \mu^Z(\cdot;ds,dy). 
\end{eqnarray*} \normalsize
It is possible to calculate the integral with respect to the compensator $\nu(\cdot;ds,dy)=N(dy) \ ds$ instead of the measure itself `under the expectation', since the integrands are of class $F_p^2$ of Ikeda-Watanabe (\cite{ikedawat}):
\[
F_p^2=\left\{ f(s,y,\omega):f \text{ is predictable, }\Bbb{E}\int_0^t\int_{\Bbb{R}} |f(s,y,\cdot)|^2 N(dy)ds \text{ for every } t>0\right\}.
\]
One obtains this, because $1_{\{|Sy|<1\}} \cdot 1_{\{ \left|\log(1+(\lambda/\delta) \ y^2)\right|<1 \}} - 1_{\{|y| < 1\}}$ is zero near the origin and bounded and $\log(1+\frac{\lambda}{\delta} y^2) \leq (\lambda/\delta)\cdot y^2$ for $\left|(\lambda/\delta)\cdot y^2\right|<1$. \newline
For $t$ tending to zero (and multiplying with $-1$) we obtain by using Lebesgue's theorem again twice\footnotesize
\begin{eqnarray*} 
&&p\left( \binom{g}{v} ,\binom{\xi_1}{\xi_2} \right) =\\
&&\hspace{9mm}-i \xi_1 \left(\ell S + S \int_{\Bbb{R} \backslash \{ 0\}} y \cdot (1_{\{|Sy|<1\}} \cdot 1_{\{ \left|\log(1+(\lambda/\delta) \ y^2)\right| <1 \} } - 1_{\{|y| < 1\}}) \ N(dy) \right)\\
&&\hspace{9mm}- i\xi_2 \left(\frac{\beta}{S^2} + \log \delta + \int_{\Bbb{R} \backslash \{ 0\}} \log(1+\frac{\lambda}{\delta} y^2) \cdot (1_{\{|Sy|<1\}} \cdot 1_{\{ \left|\log(1+(\lambda/\delta) \ y^2)\right|<1 \}} ) \ N(dy) \right) \\
&&\hspace{9mm}+\frac{1}{2} \xi_1^2 S^2 Q \\
&&\hspace{9mm}-\int_{\Bbb{R}^2 \backslash \{0 \}} \left( e^{i(z_1,z_2)\xi}-1-i z'\xi \cdot (1_{\{|z_1|<1\}} \cdot 1_{\{ |z_2|<1 \}} ) \right) 
  \tilde{N} \left( \binom{g}{S} ,dz\right), \\
\end{eqnarray*} \normalsize
where $\tilde{N}$ is the image measure
\begin{eqnarray*}
\tilde{N} \left( \binom{g}{S} ,dz \right) = N \left(\binom{S\cdot}{\log(1+(\lambda/\delta) \ \cdot^2)} \in dz \right).
\end{eqnarray*}
And by writing the starting point as $S=\exp(v/2)$ we obtain the result. \hfill $\square$

It is an advantage of our approach that, having calculated the symbol, one can write down the (extended) generator and the semimartingale characteristics at once. For the reader's convenience we recall the definition of the extended generator (cf. Definition (7.1) of \cite{vierleute}):

\textbf{Definition:}
An operator $G$ with domain $\mathcal{D}_G$ is called \textbf{extended generator} of a Markov semimartingale $X$ if $\mathcal{D}_G$ consists of those functions $f\in\mathcal{B}(\Bbb{R}^d)$ for which there exists a function $Gf\in\mathcal{B}(\Bbb{R}^d)$ such that the process 
\begin{eqnarray*}
C_t^f:= f(X_t)-f(X_0)-\int_0^t Gf(X_s) \ ds
\end{eqnarray*}
is well defined and a local martingale.

Combining Theorem 4.4 of \cite{mydiss} and Theorem 7.16 of \cite{vierleute} we obtain:

\textbf{Corollary 1:}
The extended generator $G$ on $C_b^2(\Bbb{R}^2)$ of the process $(X^{(1)},X^{(2)})'=(G,\log(\sigma^2))'$ can be written as \footnotesize
\begin{eqnarray*}
&&Gu(x)= \\
&&\hspace{2mm} \partial_1 u(x) \left(\ell e^{x_2/2} + e^{x_2/2} \int_{\Bbb{R} \backslash \{ 0\}} y \cdot (1_{\{\left| e^{x_2/2}y\right|<1\}} \cdot 1_{\{ \left|\log(1+(\lambda/\delta) \ y^2)\right|<1 \}} - 1_{\{|y| < 1\}}) \ N(dy) \right)\\
&&\hspace{2mm}+ \partial_2 u(x) \left(\frac{\beta}{e^{x_2}} + \log \delta + \int_{\Bbb{R} \backslash \{ 0\}} \log(1+\frac{\lambda}{\delta}y^2) \cdot (1_{\{\left| e^{x_2/2}y\right|<1\}} \cdot 1_{\{ \left| \log(1+(\lambda/\delta) \ y^2)\right|<1 \}} ) \ N(dy) \right) \\
&&\hspace{2mm}+ \partial_1 \partial_1 u(x) e^{x_2} Q \\
&&\hspace{2mm}+\int_{\Bbb{R}^2 \backslash \{0 \}} \Big( u(x-y)-u(x)+ y'\nabla u(x) \cdot (1_{\{|y_1|<1\}} \cdot 1_{\{ |y_2|<1 \}} ) \Big) \tilde{N} \left( x ,dy \right) \\
\end{eqnarray*} \normalsize
with the $\tilde{N}$ from above.

Writing $D(A)$ for the domain of the generator $A$ of the process we have $D(A)\subseteq \mathcal{D}_G$ and the operators $A$ and $G$ coincide on $D(A)$.

\textbf{Corollary 2:} The semimartingale characteristics $(B,C,\nu)$ of the process $(X^{(1)},X^{(2)})'=(G,\log(\sigma^2))'$ are \footnotesize
\begin{eqnarray*}
&&\hspace*{-6mm}B_t^{(1)}=\int_0^t \left(\ell e^{\frac{X^{(2)}}{2}} + e^{\frac{X^{(2)}}{2}} \int_{\Bbb{R} \backslash \{ 0\}} y \cdot (1_{\left\{\left| e^{\frac{X^{(2)}}{2}}y\right|<1\right\}} \cdot 1_{\{ \left|\log(1+(\frac{\lambda}{\delta}) \ y^2)\right|<1 \}} - 1_{\{|y| < 1\}}) \ N(dy) \right)  ds  \\
&&\hspace*{-6mm}B_t^{(2)} =\int_0^t \left(\frac{\beta}{e^{X^{(2)}}} + \log \delta +  \int_{\Bbb{R} \backslash \{ 0\}} \hspace{-3mm}\log(1+\frac{\lambda}{\delta}y^2) \cdot (1_{\left\{\left|e^{X^{(2)}/2}y\right|<1\right\}} \cdot 1_{\{ \left|\log(1+(\frac{\lambda}{\delta}) \ y^2)\right|<1 \}} ) \ N(dy) \right)  ds  \\
&&\hspace*{-6mm}C_t = \int_0^t \left( \begin{array}{cc} e^{X^{(2)}}Q & 0 \\ 0 & 0 \end{array} \right)  ds\\
&&\hspace*{-6mm}\nu(\cdot;ds,dy)= \tilde{N}(X_s(\cdot),dy) \ ds 
\end{eqnarray*} \normalsize
with the $\tilde{N}$ from above.

Remark: A different approach to calculate the characteristics of the COGARCH process is described in \cite{kallsenves}. Furthermore our results are related to earlier work of B. Rajput and J. Rosinski. In their interesting article \cite{raj-ros} they derive under certain restrictions a representation of the characteristic function of processes of the form $X_t=\int_0^t f(t,s) \ dZ_s$ where $f$ is a deterministic function and $Z$ is a L\'evy process.

\textbf{Acknowledgments:} Most of this work was done as a part of my PhD thesis, written under the guidance of Ren\'e L. Schilling to whom I am deeply grateful. Financial Support by the DFG-SFB 823 is gratefully acknowledged. Furthermore I would like to thank an anonymous referee for carefully reading the manuscript and offering useful suggestions which helped to improve the paper.

\makealttitle

\end{document}